\documentclass[12pt,a4paper]{amsart}
\usepackage{amssymb}
\usepackage[all]{xy}
\newtheorem{thm}{Theorem}[section]
\newtheorem{prop}{Proposition}[section]
\newtheorem{lemma}{Lemma}[section]
\theoremstyle{remark}
\newtheorem{remark}[section]{Remark}
\newcommand{\R}{{\mathbb R}}\newcommand{\Z}{{\mathbb
Z}}\newcommand{\C}{{\mathbb C}}
 
 \newcommand{\cO}{{\mathcal O}}
\newcommand{\cM}{{\mathcal M}}
\let\o=\operatorname
\def\bysame{$\raise.2em\hbox to 3em{\hrulefill}$\thinspace, } 
\newcommand{\Ext}{\operatorname{Ext}}
\newcommand{\Yext}{\operatorname{Yext}}
\newcommand{\cS}{\mathcal {S}}
\newcommand{\A}{\mathcal A}
\newcommand{\RB}{\mathbf {R}}

\newcommand{\PO}{\mathcal {P}}
\newcommand{\Pic}{\operatorname{Pic}}
\newcommand{\Id}{\operatorname{Id}}
\newcommand{\al}{\alpha}
\newcommand{\E}{\mathcal {E}}
\let\lfd=\longrightarrow
\begin{document}
\title{COMPLEX LAGRANGIAN EMBEDDINGS OF\\[10pt]
MODULI SPACES OF VECTOR BUNDLES}
\author{U. Bruzzo \& \ F. Pioli}
\thanks{E-Mail addresses: {\tt bruzzo@sissa.it}, {\tt pioli@sissa.it}.}
\maketitle
\begin{center}\baselineskip=12pt
\par\vskip-4mm\par
{Scuola Internazionale Superiore di Studi Avanzati } \par
{(SISSA), Via Beirut 2-4, 34014 Trieste, Italy}
\end{center}
\par\vspace{8mm}\par
\begin{quote}\footnotesize\baselineskip=14pt
{\sc Abstract.} By means of a Fourier-Mukai transform we
embed   moduli spaces $\cM_C(r,d)$ of stable bundles on an algebraic curve
$C$ of
genus $g(C)\ge 2$
as isotropic subvarieties  of moduli spaces   of $\mu$-stable
bundles on the Jacobian variety
$J(C)$. When $g(C)=2$ this provides new examples
of special Lagrangian submanifolds.
\end{quote}
\par\addvspace{8mm}\par
\section{Introduction}
Throughout this paper we shall fix $\C$ as the ground field.  Let $C$ be a
smooth
algebraic curve of genus $g>1$, and denote by $J(C)$ its Jacobian variety
and by
$\Theta\in H^2(J(C),\Z)$ the cohomology class corresponding to the theta
divisor.
Fix coprime positive integers $r$, $d$ such that
$d>2rg$, and let $\cM_C(r,d)$ be the moduli space of stable vector  bundles
on $C$ of
Chern character $(r,d)$. We show that  $\cM_C(r,d)$ can  be embedded
as an isotropic holomorphic submanifold of the complex symplectic variety
$\cM^\mu_{J(C)}(r,d)=\cM_{J(C)}^\mu(d+r(1-g),-r\Theta,0,\dots,0)$ ---
the moduli space of $\mu$-stable
vector bundles on $J(C)$ with Chern character
$(d+r(1-g),-r\Theta,0,\dots,0)$ (cf. Theorem \ref{tiram}
for a precise statement). When $g(C)=2$ one has
$\dim\cM^\mu_{J(C)}(r,d)=2\dim\cM_C(r,d)$,
and by using the  hyper-K\"ahler structure of $\cM^\mu_{J(C)}(r,d)$,
one can choose on this space a complex structure such that $\cM_C(r,d)$
embeds as a special Lagrangian submanifold, thus providing new examples of
such objects.

We recall a few facts about the Fourier-Mukai transform in the context of
Abelian
varieties \cite{Mu1}. Let $X$ be an Abelian variety  and
$\widehat{X}=\Pic^0 (X)$ its dual variety. Let $\PO$ be the normalized
Poincar\'e
bundle on $X\times
\widehat{X}$. The Mukai functor is  defined as
\begin{gather*}
\RB \cS \colon D(X)\to D(\widehat{X})\\
\RB \cS (-) =\RB\pi_{\widehat{X},\ast } (\pi_X^\ast  (-)\otimes \PO)
\end{gather*}
where $D(X)$ and $D(\widehat{X})$ are  the bounded derived categories of
coherent sheaves on $X$ and $\widehat{X}$, respectively.
Mukai has shown that the functor $\RB \cS$
is invertibile and preserves families of sheaves (cf. \cite{Mu1,Mu3}). If $E$
is a $\text{WIT}_i$ sheaf on $X$, that is, a sheaf whose transform is
concentrated in
degree $i$, then the functor $\RB \cS$ preserves the Ext groups:
$$\Ext^j_X (E,E) \cong \Ext^j_{\widehat {X}} (\hat E,\hat E)
\quad \text{for every } j, $$
where $\hat E$ indicates the transform of $E$.

Let $C$ be a smooth projective curve of genus $g> 1$ and $J(C)$ the
Jacobian of $C$.
If we fix a base point $x_0$ on $C$,  and let $\al_{x_0} \colon C\to J(C)$
be the
Abel-Jacobi embedding given by
$\al_{x_0}(x)={\mathcal{O}}_C (x-x_0)$,
the normalized Poincar\'e bundle $\PO_C$ on
$C\times J(C)$ is the pullback of the Poincar\'e bundle on $J(C)\times
J(C)$, where
we identify $J(C)$ with $\widehat{J(C)}$ via the isomorphism $-\phi_\Theta
\colon
J(C) \to
\widehat{J(C)}$. The Poincar\'e bundle on $C\times J(C)$ gives rise to a
derived
functor (which is not invertible):
\begin{gather*}
\RB \Phi_C\colon D(C)\to D(J(C))\\
\RB \Phi_C (-) =\RB\pi_{J(C),\ast } (\pi_C^\ast  (-)\otimes \PO_C)\,.
\end{gather*}
 Since  $\al_{x_0}$ is a closed immersion we have a natural isomorphism of
functors
\begin{equation}\label{relafond}
\RB \Phi_C \cong \RB \cS \circ \al_{x_0,\ast }. \end{equation}
Thus the study of the transforms of bundles $F$ on $C$ with respect  to
$\RB \Phi_C$
is equivalent to  studying the transforms of sheaves of pure dimension $1$
of the form
$\al_{x_0,\ast } (F)$ with respect to $\RB \cS$. We recall the following
fact which is
proven in \cite{Li}.
\begin{prop} If $E$ is a stable bundle on $C$ of rank $r$ and degree $d$
such that
$d>2rg$, then  $E$ is {\rm WIT}$_0$, and
the transformed sheaf $\hat E = \RB^0 \Phi_C (E)$ is locally free and
$\mu$-stable with respect to the theta divisor on $J(C)$.
\end{prop}
\par\addvspace{8mm}\par
\section{Complex Lagrangian embeddings}
If we consider the moduli
space $\cM_C(r,d)$  of stable bundles  of rank $r$ and degree $d$  on a
projective
smooth curve of genus $g>1$ such that $d>2rg$ and $r,d$ are coprime,  the
functor
$\RB \Phi_C$ gives rise to an injective morphism
$$\tilde \jmath \colon \cM_C(r,d) \to \cM^\mu_{J(C)}(r,d)=
\cM^{\mu}_{J(C)}(d+r(1-g),-r\Theta,0,\dots,0)$$
where the sheaves in $\cM_{J(C)}^\mu(r,d)$ are stable with respect to the
polarization $\Theta$.

Before studying the morphism $\tilde \jmath$ we need to recall some  elementary
facts about the Yoneda product of Ext groups. Let $\A$ be an abelian
category with
enough injectives. The elements of $\Ext^1_\A (E,E)$ are identified with
equivalence
classes of exact sequences $0\to E \to F\to E \to 0$ with respect to the usual
relation. This can be generalized to the groups $\Ext^2_\A (E,E)$ as
follows. We
refer to
\cite{HiSt} for proofs and details.

Consider the following commutative diagram with exact rows:
\begin{equation}
\xy
\xymatrix
{
E:\quad 0 \ar[r] &B \ar[d]^{\Id_B} \ar[r]& G_1 \ar[d]
\ar[r] &G_2 \ar[d]\ar[r] & A \ar[d]^{\Id_A}\ar[r]& 0 \\ E^\prime :\quad 0
\ar[r]& B
\ar[r]& G^\prime_1 \ar[r]& G^\prime_2 \ar[r] & A\ar[r] &0. }
\endxy
\label{triext}
\end{equation}
We write $E \twoheadrightarrow E^\prime$ when such a diagram holds.  The
relation
$\twoheadrightarrow$ is not symmetric, but it generates the following
equivalence
relation:
$E\sim E^\prime$ if and only if there exists a chain of sequences $E_0,
E_1, \dots,
E_k$ such that
$$E=E_0 \twoheadrightarrow E_1 \twoheadleftarrow E_2 \twoheadrightarrow \dots
\twoheadleftarrow E_k = E^\prime.$$ Let $\Yext_\A^2 (-,-)$ the set of such
equivalence classes.

There is an isomorphism
$$ \Yext_\A^2 (-,-)\cong \Ext^2_\A (-,-).$$
{}From now on we shall identify the above groups. Observe that  the identity of
$\Ext_\A^2 (A,B)$ is given by the class of the sequence
$$0\lfd B\stackrel{\Id_B}\lfd B
\stackrel{0}\lfd A \stackrel{\Id_A}\lfd A \lfd 0 .$$
Moreover the Yoneda product
$$\Ext^1_\A (B,A) \times \Ext^1_\A (A,C) \to \Ext^2_\A(B,C) $$
is obtained in the following way: let $E$ and $E^\prime $ be two elements of
$\Ext^1_\A (B,A)$ and $\Ext^1_\A (A,C)$ represented respectively by the
sequences
\begin{equation*}
E: \quad 0\lfd A\stackrel{\nu}\lfd F \stackrel{p}\lfd B \lfd 0 \end{equation*}
\begin{equation*}
E^\prime:\quad 0\lfd C\stackrel{i}\lfd G \stackrel{\lambda}\lfd A \lfd 0.
\end{equation*}
Then the class of the exact sequence
\begin{equation*}
0\lfd C\stackrel{i}\lfd G \stackrel{\nu\circ \lambda}\lfd F \stackrel{p}\lfd B
\lfd 0
\end{equation*}
in $\Ext_\A^2 (B,C)$ is the image of $E$, $E^\prime$ with respect to the
Yoneda
product.

We shall also need to introduce a moduli space of stable sheaves in
Simpson's sense
\cite{Simp}. For simplicity we denote the Abel-Jacobi map as $j\colon
C\to J(C)$. Observe that if $E$ is a stable bundle on $C$ then $j_\ast(E)$
is a
stable sheaf  of pure dimension 1 on $J(C)$ with respect to the polarization
$\Theta$.  Let
${\mathcal M}_{J(C)}^{\hbox{\tiny pure}}(r,d)$ be the moduli space of all
stable
pure sheaves on
$J(C)$ with   Chern character $(0,\dots,0,r\Theta,d+r(1-g))$. If $\E$ is a
flat
family of vector bundles on $C$ parametrized by a Noetherian scheme $S$, then
$j_{S,\ast}(\E)$ is a  flat family of sheaves on $J(C)\times S$ over
$S$, where $j_S\colon C\times S\to J(C)\times S $ is the embedding $j\times
\Id_S$.
Therefore one has a morphism of moduli spaces
\begin{equation}j_\ast : \cM (r,d) \to \cM_{J(C)}^{\hbox{\tiny pure}}(r,d)\,.
\label{am}\end{equation}
\begin{lemma} The morphism $\tilde \jmath \colon \cM_C(r,d)
\to M^{\mu}_{J(C)}(r,d)$ is an immersion (i.e., its tangent map
is injective). \end{lemma}
\begin{proof}
From the isomorphism given by Eq.~(\ref{relafond}) and recalling that the
transform
$\RB\cS$ preserves the $\Ext$ groups of WIT sheaves, it is enough to show
that the
same claim holds for the morphism (\ref{am}).
By the very definition of the Kodaira-Spencer map, the tangent map to $j_\ast$
may be identified with the
map $$\Ext^1_{C} (E,E)\stackrel{\phi}\hookrightarrow \Ext^1_{J(C)}
(j_\ast(E),j_\ast (E))$$
 obtained in the following way. Let
\begin{equation}\label{seqa}
A: \quad 0\lfd E \lfd F \lfd E \lfd 0
\end{equation}
be a sequence representing an element of $\Ext^1_{C} (E,E)$. If we apply the
functor $j_\ast $ to the above sequence we obtain the exact sequence
\begin{equation}\label{seqb}
B: \quad 0\lfd j_\ast (E) \lfd j_\ast (F) \lfd j_\ast (E) \lfd 0.
\end{equation}
One checks immediately that the map $\phi ([A]) = [B]$ is well defined. If
$\phi
([A]) = 0$ then $\phi ([A])$ is represented by the extension
\begin{equation}\label{spez1}
\quad 0\lfd j_\ast  (E) \lfd j_\ast (E)\oplus j_\ast (E) \lfd j_\ast  (E)
\lfd 0. \end{equation}
Now applying the functor $j^\ast $ to the above sequence and noting
that $j^\ast (j_\ast (H))$ $\cong H$ for every vector bundle $H$ on $C$ we
obtain the split
exact sequence
\begin{equation}\label{spez2}
\quad 0\lfd E \lfd E\oplus E\lfd E \lfd 0. \end{equation}
Therefore $\phi ([A]) = 0$ implies $[A]=0$ and $\phi$ is injective. \end{proof}
Mukai proved that the moduli space of simple sheaves on an abelian surface
$X$ is
symplectic; more precisely, the Yoneda pairing
$$\upsilon\colon\Ext_X^1(E,E)\times \Ext_X^1(E,E) \to
\Ext_X^2(E,E)\cong \mathbb C$$
defines a  holomorphic symplectic form on
the moduli of simple sheaves on  $X$ (cf. \cite{Mu2,Mu4}).
When $\dim X=2n>2$ to define a symplectic form on the smooth locus of the
moduli
space  one needs to choose a symplectic form $\omega$ on $X$. The
symplectic form   is then defined by the compositions (cf.~\cite{K})
\begin{eqnarray}
\Ext_X^1(E,E)\times\Ext_X^1(E,E) &\to& \Ext_X^2(E,E)
\stackrel{\o{tr}}{\relbar\joinrel\relbar\joinrel\to}
H^2(X,\cO_X) \nonumber \\  &\stackrel{\sim}{\to}&
H^{0,2}(X,\C)
\stackrel{\lambda}{\relbar\joinrel\relbar\joinrel\to}
H^{n,n}(X,\C) \cong\C \label{sf}\end{eqnarray}
where $\o{tr}$ is the trace morphisms and  the map $\lambda$ is obtained by
wedging
by $\omega^{n}\wedge\bar\omega^{n-1}$.
\begin{thm}\label{tiram} If $g(C)$ is even, and the map $\tilde{\jmath}$
embeds
$\cM(r,d)$ into the smooth locus $\cM_{J(C)}^0(r,d)$
of $\cM_{J(C)}^\mu(r,d)$, then the subvarieties $\cM_C(r,d)$ are isotropic
with respect to any of the symplectic forms defined by equation (\ref{sf}).
In particular, when $g(C)=2$ the subvarieties $\cM_C(r,d)$ are Lagrangian
with respect to the Mukai
form of $\cM^{\mu}_{J(C)}(r,d)$. \end{thm}
\begin{proof}
Since $\cM_{J(C)}^0(r,d)$ is smooth, and
$\tilde\jmath\colon \cM(r,d)\to  \cM_{J(C)}^0(r,d)$ is injective and is an
immersion, it is
also an embedding.
Now, let $E\in \cM_C(r,d)$. It is
enough to show that the Yoneda product
\begin{eqnarray*}
\Ext^1_{J(C)} (j_\ast  (E), j_\ast  (E))
&\times& \Ext^1_{J(C)} (j_\ast  (E), j_\ast  (E)) \\ && \lfd \Ext^2_{J(C)}
(j_\ast  (E), j_\ast  (E)) \end{eqnarray*}
vanishes when applied to pairs $([A],[B])$ of    elements  in
$\Ext^1_{J(C)} (j_\ast  (E),$
$j_\ast  (E))$ where $[A]$ and $[B]$ are represented,  respectively, by
the sequences
\begin{equation*}
A:\qquad 0\lfd j_\ast (E) \stackrel{\nu}\lfd j_\ast (F) \stackrel{ p}\lfd
j_\ast (E)
\lfd 0 \end{equation*}
\begin{equation*}
B: \qquad 0\lfd j_\ast (E) \stackrel{i}\lfd j_\ast (G) \stackrel{
\lambda}\lfd j_\ast (E)
\lfd 0 \end{equation*}
with $F,G\in \cM_C(r,d)$. It is enough to remark that the product of the
classes of the
sequences of sheaves on $C$
$$ 0 \to E \to F\to E \to 0\,,\qquad 0 \to E \to G \to E \to 0$$
is  zero for dimensional reasons, and apply the functor $j_\ast$.

In the case $g(C)=2$ the moduli space is smooth by the  results
in \cite{Mu2}; moreover,
$$\dim\cM_{J(C)}^\mu(r,d)=2(r^2+1)=2\dim\cM_C(r,d)\,.$$
\end{proof}
\begin{remark}
If we consider the moduli space $\cM_C(r,\xi)$ of stable bundles
on $C$ of rank $r$ and fixed determinant isomorphic to $\xi$, then
the result is trivial: the variety $\cM_C(r,\xi)$ is Fano, so that
it carries no holomorphic  forms.
\hfill $\blacktriangle$\end{remark}

\par\addvspace{8mm}\par
\section{The case $g(C)=2$\label{exb}}
In this section we elaborate on the case $g(C)=2$.
One can characterize situations where the moduli space
$\cM_{J(C)}^\mu(r,d)$ is compact. This happens for instance
in the following case.
\begin{prop} Assume $g(C)=2$, $d>4r$ and that $\rho=d-r$ is a prime number.
Then every
Gieseker-semistable sheaf on $J(C)$ with Chern character $(d-r,-r\Theta,0)$
is $\mu$-stable. Moreover, if $d>r^2+r$, every such sheaf is locally free
(this always happens when $r=1,2,3$).
\end{prop}
\begin{proof}
Since $d-r$ is prime, every sheaf in $\cM_{J(C)}(r,d)$ is properly
stable.
Let $[F]\in\cM_{J(C)}(r,d)$ and assume that the subsheaf
$G$ destabilizes $F$. Let $\o{ch}(G)=(\sigma,\xi,s)$. Standard
computations show that if $F$ is not $\mu$-stable then
$$\frac{\xi\cdot\Theta}{\sigma}=-\frac{2r}{\rho}\qquad\text{and}\qquad s<0\,.$$
Setting $n=\xi\cdot\Theta$ we have $\vert n\vert =2r\sigma/\rho$, with
$\sigma<\rho$
and $\rho>3r$.
This is impossible whenever $\rho$ is prime.

The statement about local freeness follows from the Bogomolov inequality.
\end{proof}

In the case $g(C)=2$ the complex Lagrangian embedding
$\tilde \jmath\colon \cM_C(r,d)$ $\to \cM^{\mu}_{J(C)}(r,d)$
provides new examples of \emph{special Lagrangian submanifolds.}
We refer to \cite{HL,ML} for the definition and the main properties of these
objects. Now, if $X$ is a hyper-K\"ahler manifold of complex dimension $2n$,
let $I$, $J$, $K$ be three
complex structures compatible with the hyper-K\"ahler metric, and such that
$IJ=K$. Let $\omega_I$,  $\omega_J$, $\omega_K$ be the corresponding
K\"ahler forms. Then the 2-form $\Omega=\omega_I+i\omega_J$ is a
holomorphic symplectic
form in the complex structure $K$. It is easy to check that a $K$-complex
$n$-dimensional submanifold
which is  Lagrangian with respect to $\Omega$ is special Lagrangian in the
structure $J$
\cite{H2}.

One should notice that via the
Hitchin-Kobayashi
correspondence (which identifies $\mu$-stable bundles on a K\"ahler
manifold with irreducible
Einstein-Hermite bundles, cf.~\cite{K}), the space $\cM^{\mu}_{J(C)}(r,d)$
acquires a
hyper-K\"ahler structure, compatible with a K\"ahler form provided by the
Weil-Petersson metric, and
with a holomorphic symplectic form which may be identified with the Mukai
form \cite{I}.

Therefore we
obtain the following result.
\begin{prop} The space $\cM^{\mu}_{J(C)}(r,d)$ has a complex structure
such that $\tilde \jmath\colon \cM_C(r,d)$ $\to \cM^{\mu}_{J(C)}(r,d)$ is a
special Lagrangian
submanifold.
\end{prop}
The elements of the Jacobian variety $J(C)$ act on the embedding $j\colon
C\to J(C)$
by translation, so that for every $x\in J(C)$ we have a special Lagrangian
submanifold
$\tilde \jmath_x\colon \cM_C(r,d)\to \cM^{\mu}_{J(C)}(r,d)$. This provides
a family of deformations of $\tilde \jmath(\cM_C(r,d))$ through special
Lagrangian
submanifolds. As one easily shows,  this  embeds $J(C)$ into the moduli space
$\cM_{SL}$
of special Lagrangian deformations of $\tilde \jmath(\cM_C(r,d))$
(notice that $\dim_{\R}\cM_{SL}=b_1(\cM_C(r,d))=4=\dim_{\R}(J(C))$) \cite{N}.
The case   $r=1$ is somehow trivial because
$\cM^{\mu}_{J(C)}(1,d)\simeq J(C)\times J(C)$ by
a result of Mukai \cite{Mu1}.

\medskip{\bf Acknowledgements.} We thank A.~Maciocia and M.S.~Narasimhan
for useful
suggestions or remarks.


\newpage
\begin{thebibliography}{99}
\frenchspacing

\bibitem{HL} Harvey, R., and Lawson Jr., H.B., {\it Calibrated geometries,}
Acta
Math. {\bf 148} (1982), 47--157.

\bibitem{HiSt} Hilton P. J., and Stammbach U., \emph{A Course in
Homological Algebra},
Grad. Texts Math. {\bf 4}, Springer-Verlag, Berlin-Heidelberg-New York (1971).

\bibitem{H2} Hitchin, N., {\it The moduli space of complex Lagrangian
submanifolds,}
{\tt math/9901069.}

\bibitem{I} Itoh, M., {\it Quaternion structure of the moduli space of
Yang-Mills
connections,} Math. Ann. {\bf 276} (1987), 581--593.

\bibitem{K} Kobayashi, S., {\it Differential Geometry of Complex Vector
Bundles,}
Princeton University Press, Princeton 1987.

\bibitem{Li} Li, Y., \emph{Spectral curves, theta divisors and Picard bundles},
Int. J. Math. {\bf 2} (1991), 525-550.

\bibitem{ML} McLean, R., {\it Deformations of calibrated submanifolds,}
Comm. Anal. Geom. {\bf 6} (1998), 705--747.

\bibitem{Mu1} Mukai, S., {\it Duality between $D(X)$ and $D(\hat X)$  with its
application to Picard sheaves,} Nagoya Math. J. {\bf 81} (1981), 153--175.

\bibitem{Mu2}\bysame \emph{Symplectic structure of the moduli space  of
sheaves on
an abelian or $K3$ surface}, Inv. Math. {\bf 77} (1984), 101-116.

\bibitem{Mu3}\bysame \emph{Fourier functor and its applications to the
moduli of
bundles on an abelian variety}, Adv. Pure Math. {\bf 10} (1987), 515-550.

\bibitem{Mu4}\bysame \emph{Moduli of vector bundles on $K3$ surfaces, and
symplectic manifolds}, Sugaku Expositions, {\bf 1}, No.~2 (1988), 139-174.

\bibitem{N} Newstead, P.E., {\it Topological properties of some spaces of
stable
bundles,} Topology {\bf 6} (1967), 241--262.

\bibitem{Simp} Simpson, C.T., {\it Moduli of representations of the fundamental
group of a smooth projective variety I,} Publ. Math. IHES {\bf 79} (1994),
47--129.

\end{thebibliography}
\end{document}